%% file: limit22.12.09.tex
\documentclass[12pt]{article}
\usepackage{geometry}                % See geometry.pdf to learn the layout options. There are lots.
\usepackage{amsmath}
\usepackage{amscd}
\usepackage{amsfonts}
\usepackage[english]{babel}
\usepackage[latin1]{inputenc}
\usepackage[pdftex]{graphicx}
\usepackage{amsthm, amssymb}
\usepackage{setspace}

\usepackage[all]{xy}
\usepackage{figlatex}

\DeclareGraphicsExtensions{.pdf} 
\usepackage{babel}
\usepackage{verbatim}

\newcommand{\BZ}{{\mathbb{Z}}}
\newcommand{\BN}{{\mathbb{N}}}
\newcommand{\BR}{{\mathbb{R}}}

\newcommand{\BQ}{{\mathbb{Q}}}

\newcommand{\BG}{{\mathbb{G}}}

\newcommand{\OO}{{\mathcal{O}}}

\newcommand{\gb}{\beta}

\newcommand{\gC}{\Gamma}
\newcommand{\gc}{\gamma}

\newcommand{\gS}{\Sigma}

\newcommand{\gep}{\epsilon}

\newcommand{\ga}{\alpha}
\newcommand{\gt}{\tau}

\newcommand{\ti}[1]{\tilde{#1}}

\newcommand{\Ad}{\text{Ad}}

\newcommand{\Hom}{\text{Hom}}
\newcommand{\SL}{\mathrm{SL}}

\newcommand{\PSL}{\mathrm{PSL}}

\newtheorem{prop}{Proposition}[section]
\newtheorem{thm}[prop]{Theorem}
\newtheorem{lem}[prop]{Lemma}
\newtheorem{cor}[prop]{Corollary}
\newtheorem{conj}[prop]{Conjecture}
\theoremstyle{definition}

\newtheorem{defn}[prop]{Definition}
\newtheorem{rem}[prop]{Remark}

\newtheorem{ob}[prop]{Observation}

\title{Compactifications and algebraic completions of Limit groups}
\author{J. Barlev, T. Gelander\thanks{T.G. acknowledges the financial support from the European
Research Council (ERC)/ grant agreement 203418 and the Israeli Science Foundation.}}

\begin{document}

\maketitle

\section{Introduction}
Little is known about the general theory of dense subgroups. Yet, given a dense or nondiscrete embedding of an abstract group $\gC$ in some topological group $G$, one often deduces interesting information about $\gC$ or $G$ (see \cite{zuk,primitive,Margulis,sul,AG,dense,toti,deformations}). There are several standard methods (probabilistic, algebraic and dynamical) to produce dense free subgroups in a given $G$ (see \cite{DPSS,dense,toti}). It was shown in \cite{BGSS} that a locally compact group which contains a dense free subgroup of rank two contains a dense copy of every torsion free surface group of even genus. A natural class of groups which contains free and surface groups at the bottom levels and which can be built up from them step by step using geometric constructions %(such as the Makanin-Razborov diagram) 
is the class of limit groups, defined by Z. Sela \cite{Sel} who proved that they share the same universal theory as the free groups. It seems likely that these groups also share the same topological group completions as the free groups. In particular:

\begin{conj}\label{conj}
A locally compact group $G$ which has a dense $F_2$ admits a dense copy of every nonabelian limit group.
\end{conj}

We do not prove this conjecture in general, but we do prove it in two natural cases, the {\it compact} and the {\it algebraic}. More precisely, we show that for every nonabelian limit group $L$ there is a natural number $d(L)\ge 2$ (defined at the end of Section \ref{sec:LG}) such that the following hold:

\begin{thm}\label{thm:compact} 
Suppose that $L$ is a nonabelian limit group and $G$ a compact group containing a free group $F$ of rank $ r\leq d(L)$. Then $G$ contains a subgroup isomorphic to $L$ whose closure contains $F$. 
\end{thm}

\begin{cor} If $G$ is a compact group containing a (dense) free subgroup of rank $2$ and $L$ is any nonabelian limit group, then $G$ contains a (dense) subgroup isomorphic to $L$. 
\end{cor}

\begin{thm}\label{thm:Q_p}
Let $k$ be a local field of characteristic $0$ and let $\BG$ be an algebraic group defined over $k$. If $\BG(k)$ admits a dense free subgroup of rank $r$ then it contains a dense copy of any limit group $L$ with $d(L)\ge r$.
\end{thm}

Since semisimple groups over local fields are topologically two generated (see \cite{kuranishi} or \cite[Theorem 2.4]{dense} for $k=\BR$, or \cite{redundant} for general $k$), we deduce:

\begin{cor}\label{cor:1.5}
Let $\BG$ be a Zariski connected semisimple algebraic group over a local field $k$ of characteristic $0$. Then $\BG(k)$ admits a dense copy of every nonabelian limit group.
\end{cor}

The case where $k=\BR$ and $G$ is connected is simpler and was proved in \cite{BGSS}.

It follows in particular that every limit group admits a faithful $2$-dimensional representation over any local field $k$; take $G=\SL_2(k)$ in Corollary \ref{cor:1.5}. 

As another application we obtain an elementary proof to the fact that limit groups are primitive (this could also be deduced from \cite[Theorem 1.9]{primitive} or from \cite[Theorem 1.15]{primitive} and \cite[Theorem 0.3]{dahmani}). Indeed, Abert and Glasner \cite{AG} observed that any dense subgroup of $\PSL_2(\BQ_p)$ is primitive, since its intersection with $\PSL_2(\BZ_p)$ is maximal and has trivial core: $\PSL_2(\BZ_p)$ is maximal in $\PSL_2(\BQ_p)$ (cf. \cite[Proposition 3.14]{Pla-Rap} and \cite{Prasad}) and has trivial core by simplicity of $\PSL_2(\BQ_p)$, and since $\PSL_2(\BZ_p)$ is open, these properties are bequeathed to intersections with dense subgroups.

Softening the density assumption to nondiscreteness, we prove the analog of Conjecture \ref{conj}:
 
\begin{thm}\label{thm:nondiscrete}
Let $G$ be a locally compact group, and suppose that $G$ admits a nondiscrete free subgroup. Then $G$ admits a nondiscrete copy of any nonabelian limit group.
\end{thm}

\medskip

Next we come back to the case of surface groups, dealt with in \cite{BGSS}. Unfortunately, the proof of the general case in \cite{BGSS} has a mistake. We dedicate Section \ref{sec:surface} to correcting this mistake and to extending the result to surfaces of odd genus:

\begin{thm} \label{thm:surface}
Let $G$ be a locally compact group which contains a nondiscrete free group $F_r$ of rank $r>1$. Then for any $g\ge 2r$, $G$ has a subgroup containing $F$ isomorphic to a surface group $\gS_g$ of genus $g$. In particular, if $G$ admits a dense copy of $F_r$, then it also admits a dense copy of $\gS_{g}$ for every $g\ge 2r$.
\end{thm}

Note that $d(\gS_g)=g$, so Theorem \ref{thm:surface} does not yield Theorems \ref{thm:compact}, \ref{thm:Q_p} and \ref{thm:nondiscrete} for surface groups.

Our approach is similar to the one in \cite{BGSS}.

\medskip

\noindent {\bf Acknowledgment:}
{\it We would like to express exceptional gratitude to the anonymous referee whose corrections, remarks and suggestions lead to improvements in several parts of the paper. In particular, for exposing a gap in Section 7 of the first manuscript which required a nontrivial consideration. We feel lucky that our paper received such  dedicated and professional treatment.}

\tableofcontents

%%%%%%%%%%%%%%%%%%%%%%%%%%%%%%%%%%%%%%%%%%%%%%%%%%

\section{Limit groups}\label{sec:LG}
Limit groups were first introduced by Zlil Sela \cite{Sel} in his
solution of the Tarski problem which asserts that any two
nonabelian finitely generated (hereafter fg) free groups have the same elementary theory. Limit
groups are exactly the groups which have the same universal theory
as fg free groups (i.e.\ satisfy the same universal sentences;
sentences of the form $\forall x_1,\ldots,x_n,~\phi(x_1,\ldots ,x_n)
$ where $\phi $ is quantifier-less).

Limit groups may be characterized in several different ways, besides by their universal theory. One way, to which limit
groups owe their name and Sela's original definition (which we will not
give here, see \cite[Section 1]{Sel}), involves Gromov--Hausdorff limits of group
actions on trees. Instead we give another definition in the flavor of
limits given in \cite{BF}:

\begin{defn} \label{limit group}
Let $G$ be a fg group and $F$ a free group. A sequence $\{
f_i\}\subseteq \Hom (G,F)$ is called \textit{stable} if 
$$
 \limsup_i
 (\ker f_i) = \liminf_i(\ker f_i),
$$ 
i.e.\ if every element $g\in G$ is eventually always or eventually never in the kernel. A fg
group $L$ is a \textit{limit group} if there exists a fg group
$G$ and a stable sequence $\{ f_i\}\subseteq \Hom (G,F)$ such that
$$
 L=G/\lim_i \ker f_i.
$$

\end{defn}

Yet another definition is given in \cite{GC} where the set of
$n$-generated groups are topologized in such a way so that limit
groups are the closure of free groups in this topology.

\subsection{Residual freeness}

It turns out that limit groups have long been studied as
\textit{fg fully residually free groups}:

\begin{defn}
A {\it residually free} group is a group $G$ such that for every $g\in G\setminus\{1\}$ there exists a free group $F$ and a homomorphism $f:G\rightarrow F$ so $f(g)\neq 1$.
\end{defn}

\begin{defn}
A {\it fully residually free} group is a group $G$ such that for
every $g_1,\ldots ,g_n \in G$ there exists a free group $F$ and a
homomorphism $f:G\rightarrow F$ so $f(g_i)\neq 1$ for $1\leq i\leq
n$.\
\end{defn}

\begin{rem} Residually free does not imply fully residual free. For example $F_2\times F_2$ is clearly  residually free, however if $w,w'\in F_2\setminus\{1\}$ then $(w,1),(w',1),([w,w'],1),(1,w)$ are not separated by any homomorphism into a free group. Indeed, as $(1,w)$ commutes with $(w,1),(w',1)$ any homomorphism must map all three to some cyclic subgroup, in which case $([w,w'],1)$ is mapped to the identity.
\end{rem}

\begin{thm}{(\cite{Sel}, 4.6)}\label{frf}
A fg group $L$ is a limit group if and only if it is fully
residually free.
\end{thm}

\begin{cor} \label{subgroup}
A fg subgroup of a limit group is a limit group.
\end{cor}

%%%%%%%%%%%%%%%%%%%%%%%%%%%%%%%%

\subsection{Constructibility}

The last approach, and the one most useful to us, presents limit
groups inductively as certain graphs of simpler limit groups. Several
such constructions are known, and they all rely on certain finiteness
properties:

\begin{thm}{(\cite{Sel}, 4.4)}
Limit groups are finitely presented.
\end{thm}

\begin{thm}{(Finite Length \cite{Sel}, 5.1)}\label{finite length}
Every sequence of epimorphisms between limit groups
$$L_1\twoheadrightarrow L_2\twoheadrightarrow L_3\twoheadrightarrow \cdots $$ stabilizes, i.e.\ all homomorphisms are isomorphisms for $n$ large enough.
\end{thm}

We now present two constructions from \cite{GC}:

\begin{defn} {(\cite{GC}, 4.10)} \label{simple graph of limit groups}
A \textit{simple graph of limit groups} over a limit group $L$, is a
group $G$ which is the fundamental group of a graph of groups such
that:

\begin{enumerate}
\item Each vertex is fg.
\item Edge homomorphisms are injective and each edge group is a nontrivial abelian group whose images under both edge morphisms are maximal abelian
subgroups of the corresponding vertex groups.
\item $G$ is commutative transitive, i.e. $\forall g,h,k\in G\setminus 1, [g,h],[h,k]=1$ implies $[g,k]=1$.
\item There is an epimorphism $\phi :G\twoheadrightarrow L$ which is
    injective when restricted to each vertex group.
\end{enumerate}
\end{defn}

\begin{thm} \label{GC4.15}(\cite{GC}, 4.15)
Every freely indecomposable limit group is a simple graph of limit
groups over a proper quotient (i.e. not injective).
\end{thm}

Next we recall the definition of an HNN extension which appears in \ref{generalized double}:

\begin{defn}\label{HNN}
Let $G$ be a group, $C$ a subgroup of $G$ and $\alpha:C \rightarrow G$ a
homomorphism. The \textit{HNN extension of $G$ over $\alpha$} is the
group
$$G*_C:=\langle G,t ~| ~tc=\alpha (c)t\rangle.$$
The element $t$ is called the stable element.
\end{defn}

The following definition is a special case of \ref{simple graph of limit
groups}:

\begin{defn} (\cite{GC}, 4.4)  \label{generalized double}
A \textit{generalized double} over a limit group $L$ is a group $A*_C B$ (or $A*_C$) such that:
\begin{enumerate}

\item Both vertex groups are fg.
\item Edge homomorphisms are injective and $C$ is a nontrivial abelian group whose image under both embeddings are maximal abelian in the vertex groups.
\item There is an epimorphism $\phi :G \twoheadrightarrow L$
which is injective when restricted to the vertices.
\end{enumerate}

\end{defn}

The following theorem is the counterpart of \ref{GC4.15}:

\begin{thm}\label{GC4.6}{(\cite{GC}, 4.6 and 4.13)}
Every freely indecomposable limit group is a generalized double over
a proper quotient.
\end{thm}

The next observation shows that any non-free limit group is a generalized double over a proper quotient:

\begin{ob} \label{slide} If
$$\xymatrix{L_1 *_A L_2 \ar[r] & L'}~~\text{or}~~\xymatrix{L_1 *_A  \ar[r] & L'}$$
is a generalized double then
$$\xymatrix{ L_1 *_A (L_2 *K) \cong (L_1 *_A L_2)  *K \ar@{->>}[r] & L'*K}$$
or respectively 
$$\xymatrix{ (L_1 *K)*_A  \cong (L_1 *_A )  *K \ar@{->>}[r] & L'*K}$$
(with the obvious map) is also generalized double
with vertices $L_1$ and $L_2 *K$ in the amalgamated product case, or respectively $L_1*K$ in the HNN case. As a consequence, if a limit group decomposes as a free product  $L*K$ with $L\twoheadrightarrow L'$ a generalized double quotient for $L$ then the above map $L*K\twoheadrightarrow L'*K$ is a generalized double quotient for $L*K$. 
\end{ob}

As a consequence of \ref{finite length}, \ref{GC4.6} and \ref{slide} we obtain:

\begin{cor} \label{lgdd}
Let $L$ be a limit group. There exists a finite sequence 
$$L\twoheadrightarrow L_1 \twoheadrightarrow \cdots \twoheadrightarrow L_n\twoheadrightarrow F$$
where each $L_i\twoheadrightarrow L_{i+1}$ is a generalized double and $F$  is a fg free group. 
\end{cor}

The sequence of generalized doubles should be thought of as a construction of $L$ from $F$ using generalized doubles.

\medskip

\noindent {\bf The number $d(L)$:} {\it For a limit group $L$, we define $d(L)$ to be the maximal possible rank of a free group $F$ appearing at the end of a generalized double construction sequence \ref{lgdd}.}

%%%%%%%%%%%%%%%%%%%%%%%%%%%%%%%%%%%%%%%%%%%%%%%%%%%%

\section{A generalization of a lemma of Baumslag}

In this section $F$ denotes a free group of finite rank. Recall the following classical (and very useful) result of Baumslag \cite{B}: 

\begin{lem}[Baumslag's lemma]\label{lem:Baumslag}
Let $a_0, a_1,\ldots ,a_n,~n>0$ and $z$ be elements in a free group $F$ such
that $[z,a_i]\neq 1$ for $i=1,\ldots,n$.
Then there exists an integer $N$ such that for all integers $k_0, \ldots,
k_n$ with $|k_i|\geq N$ for $i = 0,\ldots,n-1$,
$$
 a_0z^{k_0}a_1z^{k_1}a_2\cdots z^{k_{n-1}}a_nz^{k_n}\neq 1.
$$
\end{lem}

This lemma immediately generalizes to fully residually
free groups, and in particular to limit groups.

The following notion is useful:

\begin{defn}[Dehn twist]\label{dehn twist}
Let $\gC$ be an amalgamated product $\gC=A*_C B$ or HNN extension
$\gC=A*_C$ with $C$ abelian. For every $c\in C$ the \textit{Dehn twist}
associated to $c$ is the automorphism $\tau_c:\gC\rightarrow \gC$ such that
in both cases for $a\in A$, $\tau_c(a)=a$; in the case of the
amalgamated product for $b\in B$, $\tau_c(b)=c^{-1}bc$; and in the
case of the HNN extension $\gt_c(t)=c^{-1}t$ where $t$ is the stable element.
\end{defn}

As an example of a use of Lemma \ref{lem:Baumslag} let us show (the well known fact) that a double of a
limit group amalgamated over a maximal abelian subgroup,
$L*_{C=\overline{C}} \overline{L}$ (where $L\cong \overline{L}$), is a limit group.
Let $z\in C\setminus \{1\}$ and consider the Dehn twist $\tau$ along
$z$. Let $f:L*_{C=\overline{C}}
\overline{L}\rightarrow L$ be the folding map. By the Baumslag
lemma, $f \circ \tau ^n$ is eventually faithful (i.e.\ stable with trivial stable kernel; see Definition \ref{limit group}) hence the double is
fully residually free, i.e. a limit group.

Another useful consequence of Lemma \ref{lem:Baumslag} is the following.
Let $F_{n+1}$ be the free groups on $\{x_1,\ldots,x_{n+1}\}$ and let $F_n=\langle x_1,\ldots,x_n\rangle$. Pick $a\in F_n$ nontrivial and define $\gt_a \in\text{Aut}(F_{n+1})$ to be the identity on $F_n$ and $\gt_a(x_{n+1})=ax_{n+1}a^{-1}$. Let $b\in F_n$ be an element not commuting with $a$ and let $f:F_{n+1}\to F_n$ be the homomorphism whose restriction to $F_n$ is the identity and $f(x_{n+1})=b$. Then the sequence $f\circ \gt_a^m:F_{n+1}\to F_n$ is eventually faithful. This is summarized as follows:

\begin{cor}\label{cor:F_ndouble}
Let $A$ be a maximal cyclic subgroup of $F_n$. Then 
$$
 f:F_{n+1}\cong F_n*_A(A*\langle x_{n+1}\rangle)\to F_n
$$ 
is a generalized double over $F_n$ (where $f(x_{n+1})=b\notin A$) and for every nontrivial $a\in A$ the sequence $f\circ\gt_a^m$ is eventually faithful.
\end{cor}

One can formulate many useful variants of Lemma \ref{lem:Baumslag}, for instance:

\begin{cor}\label{cor:baum1}
Let $a,b$ be noncommuting elements in a free group $F$, let $l_1,\ldots,l_m\in\BZ\setminus\{0\},~m>0$ and let $w_0,w_1,\ldots,w_m$ be arbitrary elements in $F\setminus\{ 1\}$. Then for all sufficiently large $k$,
$$
 w_mb^{-k}a^{l_m}b^{k}w_{m-1}\cdots w_1b^{-k}a^{l_1}b^kw_0\ne 1.
$$
\end{cor}

\begin{proof}
Without loss of generality assume that $b$ is the generator of the cyclic group $C_F(b)$. Hence if $w_i$ commutes with $b$, it has the form $b^{t_i}$ for some $t_i\ne 0$, and we can, for every such $i<m$, replace in the above expression $a^{l_{i+1}}b^{k}w_ib^{-k}a^{l_i}$ by $(a^{l_{i+1}}b^{t_i}a^{l_i})$. This gives a, perhaps shorter, expression whose terms are $b^{\pm k}$ and elements not commuting with $b$, alternately, $w_m$ being the only possible exception, and this expression starts with $w_mb^{-k}\ldots$. By Lemma \ref{lem:Baumslag} it is nontrivial for all large $k$.
\end{proof}

We will need a more general version of Lemma \ref{lem:Baumslag} that will allow us to deal with graphs of groups with more than one edge.

\begin{lem}[A generalized Baumslag's lemma]\label{baumslag}
Let $z_{1},\ldots,z_{l},u_{1},\ldots,u_{m}\in F$ be elements satisfying 
\begin{itemize}
\item $[u_{i},z_{j}]\neq 1$ for $1\leq i\leq m,~1\leq j\leq l$, and
\item $[u_{i}z_{j}u_{i}^{-1},z_{k}]\neq1$ for $0\leq i\leq m$ (where $u_{0}=1$), $1\leq
j,k\leq l,~~ j\neq k$.
\end{itemize}
Then there is an integer $N\in\BN$ such that every
$$ 
 g=w_{n}\cdots w_{1},
$$ 
where
\begin{itemize}
\item $w_{i}=u_{j_i},~j_i\in\{0,\ldots,m\}$ for every odd $i\le n$,
\item $w_{i}=z_{k_i}^{t_i},~k_i\in\{1,\ldots, l\}$, $t_i\in\BZ$ for every even $i<n$,
\item $w_{i_0}=u_0$ for $i_0$ odd $\Rightarrow$ $n\ge 2$, and if additionally $n>i_0>1$ then $z_{k_{i_0-1}}\ne z_{k_{i_0+1}}$,
\end{itemize}
is nontrivial provided $|t_i|\ge N,~\forall$ even $i\le n$.
\end{lem}

Let $T$ be the Cayley graph of $F$ with respect to a set
of free generators. Denote by $\partial T$ the ideal boundary of $T$, i.e.\ the
topological space whose underlying set is that of rays $r_\ga$ in $T$ up to the
equivalence relation $r_\ga\sim r_\gb$ if they eventually coincide. A basis
for the topology of $\partial T$ is given by the subsets $\{r\in
\partial T: e\in r\}$ for all oriented edges $e\in T$. $\partial T$ is a
compact, Hausdorff and totally disconnected space on which $F$ acts by
homeomorphisms. Every nonidentity element $a\in F$ acts hyperbolically on
$T$ and hence has an axis, denoted $\text{axis}(a)$ with endpoints $a^+,a^-\in\partial T$.

\begin{lem}\label{commutator assumption}
$[a,b]\neq1$ iff $\text{axis}(a)\cap \text{axis}(b)$ is finite.
\end{lem}

\begin{proof}
Suppose $\text{axis}(a)\cap \text{axis}(b)$ is infinite.
Then it contains a ray. In particular for $x$ in
the intersection and some integers $n,k$ (depending on the displacement of $a$ and $b$) we have $a^{k}x=b^{n}x$ which yields 
$[a^k,b^n]\cdot x=1$. Since $F$ acts freely on $T$, $[a^k,b^n]=1$, and the result follows since $F$ is commutative transitive.

Conversely, if $\text{axis}(a)\cap \text{axis}(b)$ is finite, for $n>|\text{axis}(a)\cap \text{axis}(b)|$, $a^n$ and $b^n$ satisfy the Schottky condition (the condition of the ping-pong lemma, see \cite[Chapter II.B.]{delaHarpe} or \cite{AG}) and hence $\langle a^n,b^n\rangle\cong F_2$.
\end{proof}

Note that $|\text{axis}(a)\cap\text{axis}(b)|<\infty\iff\{ a^+,a^-\}\cap\{b^+,b^-\}=\emptyset$.

\begin{figure}[h!]
\begin{center}

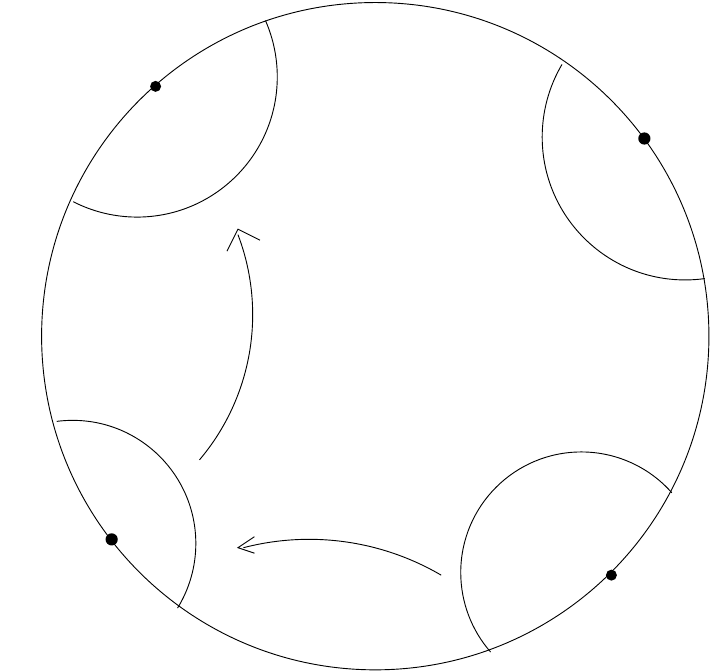

\end{center}
\end{figure}

\begin{proof}[Proof of Lemma \ref{baumslag}]
Since the case $n=1$ is trivial, assume that $n\ge 2$. We also exclude the simple special case where $n\le 3$ and $w_1=1$.
We use a ping-pong type argument. The commutator assumptions and Lemma
\ref{commutator assumption} imply that all the boundary points
$$
 \big\{z_{j}^{\pm}:~1\leq j\leq l \big\}
$$ 
are distinct, and that for every $1\leq i\leq m$, 
$$
 \big\{ u_{i}\cdot z_{j}^{\pm}:~1\leq j\leq l \big\}\cap\big\{
 z_{j}^{\pm}:~1\leq j\leq l \big\}=\emptyset. 
$$ 
Let, for each $i$,
$V_{i,j}^{\pm}$ be separating open neighborhoods of the ends $u_i\cdot z_{j}^{\pm},~i=1,\ldots,l$
respectively, which do not contain any of the $z_k^\pm$ in their closures.
Then let $V_{j}^{\pm}$ be open neighborhoods of the
$z_{j}^{\pm}$, disjoint from all the $V_{i,k}^\pm$, such that $u_{i}\cdot V_{j}^{\epsilon}\subseteq
V_{i,j}^{\epsilon}$ for $1\leq i\leq m,~\gep=\pm$ (see Figure above). Let $N\in\BN$ be such that for
$t\geq N$ we have $z_{j}^{\epsilon t}\cdot(\partial T\backslash
V_{j}^{-\epsilon})\subseteq V_{j}^{\epsilon}$ for all $j=1,\ldots,l$ (where we agree that $--$ is $+$ and $-+$ is $-$).
Assume $|t_1|,\ldots,|t_n|\geq N$.

Let us say that two points in $\partial T$ are {\it close} if they belong to one of the sets $V_j^\gep$ or $V_{i,j}^{\epsilon}$ and {\it far} if they belong to two different $V_j^\gep$. Since, for a fixed $i$ the $V_j^\epsilon,~V_{i,j}^{\epsilon},~j=1,\ldots,l,\gep=\pm$, are pairwise disjoint, two points cannot be simultaneously close and far. 

For $k=1,\ldots,n$ denote $g_{k}=w_{k}\cdots w_{1}$
and follow the orbits of, say, $z_1^+$ and $z_1^-$ under $g_1,g_2,\ldots,g_n$. One derives from the assumptions that $g_2\cdot z_1^+$ and $g_2\cdot z_1^-$ are already close and that the orbits stay close ever after, unless $w_1=1, w_2=z_1^\pm$ in which case $n\ge 4$ and the orbits become and stay close from the fourth step (or alternatively one can look at $z_2^\pm$ in this case). In particular $g_n\cdot z_1^+$ and $g_n\cdot z_1^-$ are close. 
Since $z_1^+$ and $z_1^-$ are far, this implies that $g=g_n$ is nontrivial.

\end{proof}

\begin{cor}\label{cor:bslg}
It follows from the proof that it is actually enough to require only that
$[u_{i}z_{j}u_{i}^{-1},z_{k}]\neq1$ and $[u_{i},z_{j}]\ne 1$ for $1\leq
k\leq l,\,\, j\neq k$ whenever the subword $u_{i}z_{j}$ appears in
$g$. 
\end{cor}

As before, Lemma  \ref{baumslag} is also true for limit groups.

\begin{rem}\label{bslg remark}
Our argument provides in particular an elementary proof of the original Baumslag lemma. One can obtain an even more general result by the same reasoning, but the formulation we gave is sufficient for our use in later sections.
\end{rem}

%%%%%%%%%%%%%%%%%%%%%%%%%%%%%%%%%%%%%%%%%%%%%%%%%%%%%%%

\section{A strategy for embedding limit groups}
The proofs of Theorems \ref{thm:compact} and \ref{thm:Q_p} are by induction on the length of the diagram given in Corollary \ref{lgdd}. 
Let $G$ be a group containing a subgroup $L$ and let $L'=A*_{C}B\xrightarrow{\varphi}L$ (or
$L'=A*_{C}\xrightarrow{\varphi}L$ with stable element $t$) be a generalized double over $L$. Then for every $g\in C_G(\varphi (C))$, where $C_G(\varphi(C))$ is the centralizer of $\varphi(C)$ in $G$, the rule 
$$
 a\mapsto\varphi (a),~b\mapsto g\varphi(b)g^{-1},~\text{resp.}~a\mapsto\varphi (a),~t\mapsto\varphi (t)g
$$
defines a homomorphism $h_g:L'\to G$. As a straightforward consequence of the Baire category theorem we have:

\begin{lem} \label{lem:double}
Let $G$ be a topological group and $L\le G$ a subgroup. Let
$L'=A*_{C}B\xrightarrow{\varphi}L$ (or $L'=A*_{C}\xrightarrow{\varphi}L$) be a generalized double over $L$ where $L'$ is countable. Suppose that $C_G(\varphi(C))$ admits a subset $K$ which is a Baire space with respect to the induced topology such that
\begin{itemize}
\item $L\subset\overline{\langle h_k(L')\rangle},~\forall k\in K$;
\item $\{ k\in K: h_k(l)=1\}$ is nowhere dense in $K$ for any $l\in L'\setminus\{ 1\}$. 
\end{itemize}
Then there exists an embedding
$h:L'\hookrightarrow G$ such that $h|_A=\varphi |_A$ and $\overline{h(L')}\supseteq L$.
\end{lem}

In the next two sections we will apply the following abstract proposition:

\begin{prop}\label{prop:emb-lim}
Let $G$ be a locally compact group with the property that for any dense limit subgroup $L\le G$ and a generalized double $L'=A*_{C}B\xrightarrow{\varphi}L$ (or $L'=A*_{C}\xrightarrow{\varphi}L$) over $L$, there is a subset $K\subset C_G(\varphi(C))$ such that the conditions of Lemma \ref{lem:double} are satisfied. If $G$ admits a dense free subgroup of rank $r$ then it admits a dense copy of every limit group $L$ satisfying $d(L)\ge r$.
\end{prop}  

\begin{proof}
Let $F_r$ be a dense free subgroup of rank $r$ in $G$. By the definition of  $d(L)$
there exists a generalized double construction sequence (Corollary \ref{lgdd}) for $L$ ending with a free group $F_s$ of rank $s\geq r$. 
In view of Corollary \ref{cor:F_ndouble} $F_s$ can be constructed from $F_r$ by a sequence of generalized doubles of length $s-r$. Thus $L$ admits a constructive sequence ending with $F_r$. 
Lemma \ref{lem:double} allows us to move up the
sequence inductively and show that there exists a dense embedding
of $L\hookrightarrow G$.
\end{proof}

The following will be used:

\begin{lem}\label{lem:emigroup-dense}
Let $G$ be a topological group and $c\in G$ an element generating a nondiscrete cyclic semigroup. Then $\overline{\{ c^n:n\ge N\}}=\overline{\langle c\rangle}$ for every $N\in\BN$.
\end{lem}

\begin{proof}
Denote $\mathcal{C}=\overline{\langle c\rangle}$ and 
$C_N=\{{c}^{n}:n\geq N\}$. It is enough to
show $1\in \overline{C_N}$ since the accumulation points of this set are
invariant under translations by elements from the group
$\langle{c}\rangle$. By nondiscreteness $C_N$ has some
accumulation point $b$; let ${c}^{n_{i}}\rightarrow b$. By
the above remark, ${c}^{-n_{i}}b\in \overline{C_N}$, hence
$1=\lim{}_{i}{c}^{-n_{i}}b\in C_N$.
\end{proof}

Let us now prove a simple special case of Proposition \ref{prop:extension}:

\begin{lem}\label{infusion}
Let $F$ be a free group of rank $n\ge 1$, and let $f:F\hookrightarrow G$ be an embedding of $F$ in a locally compact group. Suppose that $F$ admits an element $a$ such that the cyclic semigroup $\langle f(a)\rangle$ is nondiscrete. Then there exists an embedding
$F*\langle z\rangle\hookrightarrow G$ of the free group $F_{n+1}\cong F*\langle z\rangle$ of rank $n+1$ whose restriction to $F$ is $f$.
\end{lem}

\begin{proof}
Taking $K=\overline{\langle f(a)\rangle}$ in Lemma \ref{lem:double}, the result follows immediately from Corollary \ref{cor:F_ndouble}
and Lemma \ref{lem:emigroup-dense}.
\end{proof}

%%%%%%%%%%%%%%%%%%%%%%%%%%%%%%%%%%%%%%%%%%%%%%%%%

\section{Embedding in compact groups}\label{sec:compact}

We will now give a proof of Theorem \ref{thm:compact}.

Let $G$ be a compact group, $L\le G$ a nondiscrete limit subgroup, and $L'=A*_{C}B\xrightarrow{\varphi}L$ (or $L'=A*_{C}\xrightarrow{\varphi}L$) a generalized double over $L$. Replacing $G$ by the closure of $L$ we may assume that $L$ is dense.
By \cite{GC} 4.12, there exists $c\in C$ such that the sequence of
associated Dehn twists precomposed with $\varphi$,
$\varphi\circ\tau_c^{n}:L'\rightarrow L$ is eventually faithful. 
Let $\overline c=\varphi(c)$ and let $K$ be the closure of the cyclic group generated by $\overline c$.
%$K=\overline{\langle\overline{c}\rangle}$. 
In view of Proposition \ref{prop:emb-lim} it is sufficient to show that $K$ satisfies the conditions of Lemma \ref{lem:double}.

By the choice of $c$, for each $l\in L'\setminus\{1\}$ there is $N_l\in\BN$ such that 
$$
 n\ge N_l\Rightarrow h_{\overline{c}^n}(l)=\varphi(\gt_c^n(l))\ne 1.
$$
Hence by Lemma \ref{lem:emigroup-dense} the closed set $\{ k\in K: h_k(l)=1\}$ is nowhere dense in $K$, so the second requirement of Lemma \ref{lem:double} is satisfied.

It is left to check the first requirement of Lemma \ref{lem:double}. Let $k$ be an element of $K$.
Then $h_{k}(a)=\varphi (a),~\forall a\in A$. Therefore, it is enough to show that in the amalgamated product case, $\overline{h_{k}(B)}\supseteq\varphi(B)$, because in that case $L$ is generated by $\varphi(A)$ and $\varphi(B)$. 
Now, $h_{k}|_B$ is simply $\varphi|_B$ followed by
conjugation with $k$. Let $n_i\in\BN$ be integers so that $\overline{c}^{n_{i}}\rightarrow k$ (Lemma \ref{lem:emigroup-dense}),
and let $b\in B$ and $\overline{b}:=\varphi\left(b\right)$.
Then
$$ 
 h_{k}\left(c^{-n_{i}}bc^{n_{i}}\right)=k\varphi(c^{-n_{i}}bc^{n_{i}})k^{-1}=(k\overline{c}^{-n_{i}})
 \varphi(b)(\overline{c}^{n_{i}}k^{-1})\xrightarrow{i}
 \varphi(b),
$$ 
hence $\varphi(b)\in\overline{h_{k}\left(B\right)}$ for all $b\in B$. 
This completes the proof of Theorem \ref{thm:compact}.
\qed

\medskip

We end this section with an immediate consequence which is a special case of Theorem \ref{thm:nondiscrete}: 

\begin{cor}\label{cor:totdisc}
Let $G$ be a locally compact totally disconnected topological group containing a nondiscrete nonabelian free
subgroup $F$. Then $G$ contains a nondiscrete copy of every nonabelian limit group.
\end{cor}

\begin{proof}
Let $O$ be an open compact subgroup of $G$, and choose two noncommuting elements $a,b\in F\cap O$. Then $\langle a,b\rangle$ is a free group whose closure $H$ is compact. By Theorem \ref{thm:compact} every nonabelian limit group $L$ admits a faithful (dense) embedding in $H$.
\end{proof}

%%%%%%%%%%%%%%%%%%%%%%%%%%%%%%%%%%%%%%%%%%%%%%%%%%
\section{Embedding in algebraic groups over local fields}

In this section we explain the proof of Theorem \ref{thm:Q_p}. Let $k$ be a local field of characteristic $0$, $\BG$ a $k$ algebraic group and $G=\BG(k)$ the group of rational points. As in the previous section, given a dense limit subgroup $L\le G$, and a generalized double $L'=A*_{C}B\xrightarrow{\varphi}L$ (or $L'=A*_{C}\xrightarrow{\varphi}L$) over $L$, we have to find the Baire set $K\subset G$ which satisfies the conditions of Lemma \ref{lem:double}. In view of Proposition \ref{prop:emb-lim} this will finish the proof.

%We wish to argue along the same lines we did in the compact case.
%Thus let $L$ be a limit group and suppose that $G=\BG(k)$ admits a dense subgroup isomorphic to $L$ (we will denote it also by $L$), and let $L'=A*_CB\xrightarrow{\varphi} L$ be a generalized double over $L$. We claim that $G$ admits a dense copy of $L'$ as well.
%The argument is very similar to the one in the previous section hence we shall not repeat it, but only explain the modification needed to apply it in the context of algebraic groups over $p$-adic fields.
%Let $\BG$ be a Zariski connected $p$-adic algebraic group.  Let $\{\ga_1,\ldots,\ga_n\}$ and $\{\gb_1,\ldots,\gb_m\}$ be generating sets for $A$ and $B$ respectively. Hence $L=\langle\varphi(\ga_i),\varphi(\gb_i):i\le n,j\le m\rangle$. 

Suppose first that $k$ is nonarchimedean, i.e. a finite extension of $\BQ_p$ for some rational prime $p$. 
%By restriction of scalars we may assume that $k=\BQ_p$.
Let $\ti C$ be the Zariski identity component of the Zariski closure of $\varphi (C)$ in $G$, i.e.\ $\ti C=(\overline{\varphi(C)}^z)^\circ\cap \BG(k)$. Then $\ti C$ is a Zariski connected $p$-adic algebraic group containing $\varphi(C)$. By \cite[Lemma 3.2]{Pla-Rap} any $p$-adic open subset of $\ti C$ is Zariski dense in $\ti C$. Moreover, $G$ is a $p$-adic analytic group, hence the Frattini subgroup of some open compact subgroup of $G$ is open (see \cite[Proposition 1.14, Theorem 5.2 and Theorem 8.32]{DSMS}). It follows that, for every $m\in\BN$, the set of $m$-tuples generating dense subgroups is open in $G^m$ (see \cite[Section 5.2]{toti} for details).
As $L'$ is finitely generated, the set 
$$
 \{ f\in\text{Hom}(L',G):f(L')~\text{is dense in}~G\}
$$ 
is open in $\text{Hom}(L',G)$ (cf. \cite[Section 5]{BGSS}). (Recall that the canonical topology on $\text{Hom}(L',G)$ can be described (noncanonically) by fixing a finite generating set $\gS$ for $L'$, identifying $\text{Hom}(L',G)$ with the corresponding closed subset of $G^{|\gS|}$ and taking the induced topology.)  
%Thus the statement that $\text{Hom}(L',G)$ is open for general $L'$ follows from the case of the free group where the statement can be read as: {\it $\forall d$, the set of $d$-tuples generating a dense subgroup is open in $G^d$.} 

Let $\{a_1,\ldots,a_n\}$ and $\{b_1,\ldots,b_m\}$ be finite generating sets for $A$ and $B$ respectively ($B=\emptyset$ in the HNN case), and denote $\ga_i=\varphi(a_i)$ and $\gb_j=\varphi(b_j)$.
Set $K$ to be a sufficiently small open compact subgroup of $\ti C$ so that $\langle\varphi(\ga_i),k\varphi(\gb_i)k^{-1}:i\le n,j\le m\rangle$ is dense in $G$ for any $k\in K$. The first condition of Lemma \ref{lem:double} is then automatically satisfied. Note that by definition $K\subset \ti C\le C_G(\varphi(C))$.

By \cite{GC} 4.12 there is $c\in C$ such that $h_{\varphi(c^n)}$ is eventually faithful, and we can certainly assume that $c$ belongs to the subgroup $\varphi^{-1}(\ti C)$, the later being of finite index in $C$. 
Thus given $\gc\in\ L'\setminus\{ 1\}$ the algebraic set $\{\ti c\in \ti C:h_{\ti c}(\gc)=1\}$ is proper
and since $K$ is open and Zariski dense in $\ti C$ and $\ti C$ is $p$-adic analytic, it intersects $K$ in a nowhere dense subsets (cf. \cite[Lemma 3.2]{Pla-Rap}). Thus the second condition of Lemma \ref{lem:double} is also satisfied.

\medskip

Suppose now that $k$ is archimedean. This time, define $\ti C$ as the Hausdorff identity component of the Zariski closure of $\varphi (C)$ in $G$. Then, since $\ti C\le G^\circ$, $\overline{h_{\ti c}(L')G^\circ}=G~\forall \ti c\in\ti C$, where $G^\circ$ denotes the Hausdorff identity component of $G$. Note that by Whitney's theorem \cite[Theorem 3.6]{HR} $\ti C$ is of finite index in $\overline{\varphi (C)}^z$.
Let $\mathcal{L}'$ be the finite index subgroup $\varphi^{-1}(G^\circ)$ of $L'$.
As follows for instance from the analysis given in \cite[Section 6]{BGSS}, the set of representations of $\mathcal{L}'$ (or of any fg group) in $G^\circ$ with dense image is open in $\text{Hom}(\mathcal{L}',G^\circ)$ if $G^\circ$ is topologically perfect (see \cite[Corollary 6.2]{BGSS}), and in general it is an open set minus a countable union of nowhere dense closed subsets which are furthermore proper algebraic (see \cite[Lemma 6.4]{BGSS}) (in \cite{BGSS} it is only indicated that these sets are analytic, but their explicit description is clearly algebraic). Moreover, the map $\ti c\mapsto h_{\ti c}$ from $\ti C$ to $\text{Hom}(\mathcal{L}',G^\circ)$ is obviously analytic (and even algebraic) where the analytic and the algebraic structures on $\text{Hom}(\mathcal{L}',G^\circ)$ are, as the topological structure, canonically induced from the corresponding structures on $G$. It follows that the set 
$$
 K:=\{\ti c\in\ti C:\overline{h_{\ti c}(\mathcal{L}')}=G^\circ\}
$$ 
is also an open subset possibly minus a countable union of proper algebraic subsets. Hence $K$ 
is a Baire space. By definition $K$ satisfies the first condition of Lemma \ref{lem:double} and since it is Zariski dense in $\ti C$ it also satisfies the second condition (because, as above, we have that $h_{\ti c^n}$ is eventually faithful for some $\ti c \in\ti C$).
\qed

%\medskip

%The proof in the archimedean case above, can be extended for general Lie groups, using the algebraic structure of the semisimple group $G^\circ/\text{Rad}(G^\circ)$. Let us formulate, omitting a detailed argument:

%\begin{thm}
%Let $G$ be a real Lie group of positive dimension which contains a dense free subgroup of rank $r$. Then $G$ contains a dense copy of any limit group $L$ with $d(L)\ge r$.
%\end{thm}

%This generalizes \cite[Theorem 1.4]{BGSS}. Note that the argument given in \cite{BGSS} cannot be applied to nonconnected Lie groups.

%\begin{exam}
%Let $G=\SL_n(\BQ_p)$. Then $G$ is topologically generated by two elements. Indeed, one can show that $\SL_n(\BZ_p)$ is topologically generated by two (regular opposite) unipotents, $u_1,u_2$. Approximate $u_1,u_2$ by diagonilizable elements $x_1,x_2\in\SL_n(\BZ_p)$ which has the same eigenvalues. Then $x_1=gx_2g^{-1}$ for some $g\in\SL_n(\BQ_p)$ which we may assume not belong to $\SL_n(\BZ_p)$. If the $x_i$ are sufficiently close to the $u_i$ we have $\overline{\langle x_1,x_2\rangle}=\SL_n(\BZ_p)$ since $\SL_n(\BZ_p)$ is $p$-adic analytic. Thus the group $\overline{\langle x_1,g\rangle}$ strictly contains $\SL_n(\BZ_p)$. However, $\SL_n(\BZ_p)$ is a maximal subgroup of $\SL_n(\BQ_p)$; hence $\langle x_1,g\rangle$ is dense in $G$. It follows from \cite{BG:toti} that $G$ admits a dense $F_2$. Hence, by Theorem \ref{thm:Q_p} we obtain:

%\begin{cor}
%For every $n\ge 2$ the group $\SL_n(\BQ_p)$ admits a dense copy of any nonabelian limit group $L$.
%\end{cor}
%\end{exam}

%%%%%%%%%%%%%%%%%%%%%%

\section{The proof of Theorem \ref{thm:nondiscrete}}

%The case of Theorem \ref{thm:Q_p} where the local field $k$ is $\BR$ and $G$ is connected was proved in \cite{BGSS} using the characterization of limit groups as fully residually free groups (see Theorem 1.5 there). Arguing similarly, we prove:

We start with a simple variant of Theorem \ref{thm:Q_p}:

\begin{lem}\label{lem:S->Ad(S)}
Let $S$ be a connected semisimple Lie group and $\gC$ a nonabelian limit group. Then $S$ admits a dense copy of $\gC$ which projects faithfully to $\Ad(S)$.
\end{lem}

\begin{proof}
We argue as in the previous sections. That is, given a
dense limit subgroup $L\le S$ which projects faithfully to $\Ad(S)$, and a generalized double $L'=A*_{C}B\xrightarrow{\varphi}L$ (or $L'=A*_{C}\xrightarrow{\varphi}L$) over $L$, we show that $L'$ can be embedded densely in $S$  intersecting trivially with $Z(S)$.

First pick $c$ in $C$ such that the powers of the Dehn twist induced by $c$ precomposed by the map $L'\to L$ are eventually faithful. We then use the fact that $\Ad(S)$, being center free and semisimple, is algebraic and take $A^\circ\le \Ad(S)$ to be the Hausdorff identity component of the Zariski closure $A$ of the cyclic group generated by the image of $c$ in $\Ad(S)$.
By Whitney's theorem $A$ has finitely many connected components, hence up to replacing $c$ by some finite power  we may assume that its image in $\Ad(S)$ belongs to $A^\circ$. Since $A^\circ$ is finite dimensional abelian and connected its associated exponential map is onto, in particular, the image of $c$ lies in a $1$-parameter subgroup $\{\overline x(t)\}_t$, which can be lifted uniquely to a $1$-parameter subgroup $\{x(t)\}_t\le S$. Now for any $c'\in C$ the commutator $[x(t),c']$ is a continuous path through $1$ lying in the discrete center of $S$, hence constant, i.e. $\{x(t)\}_t\subset C_S(\varphi (C))$. This allows us to define a map $h_{x(t)}:L'\to S$ for every $t\in\BR$, as in the previous sections. For every $\gc\in L'\setminus\{ 1\}$, the map $t\mapsto \Ad(h_{x(t)}(\gc))$ is analytic, and by the choice of $c$ is not trivial for some $t$, hence nontrivial generically. Moreover since the set of dense representations is open in $\Hom(L',S)$, for all sufficiently small $t$, $h_{x(t)}(L')$ is dense in $S$, and the result follows by a Baire category argument.
 
As an alternative proof, since $S$ is connected, one can apply the argument given in \cite[proof of Theorem 1.5]{BGSS}, noting that the requirement that {\it the image of every nontrivial $\gc\in L'$ lies out of the center of $S$}, translates to an elimination of countably many additional "nowhere dense conditions".   
\end{proof}

Our next step is:

\begin{prop}\label{prop:non-sol-quo}
Let $H$ be a locally compact group which admits a compact normal subgroup $O$ such that $H/O$ is a connected nonsolvable Lie group. Let $R$ be the solvable radical of $H/O$, and $\ti R$ its pre-image in $H$. Then $H$ admits a copy of every nonabelian limit group, which projects densely and faithfully to $H/\ti R$.
\end{prop}

%\begin{lem}\label{lem:non-sol-quo}
%Let $H$ be a locally compact group which admits a compact normal subgroup $O$ such that $H/O$ is a connected nonsolvable Lie group. Then $H$ admits a nondiscrete copy of every limit group.
%\end{lem}

We will rely on several structural result for locally compact groups:

\begin{lem}\label{lem:lift-1-parameter}\label{MZ4.6} (\cite{MZ} \textit{end of 4.6}) Let $H$ be a locally compact group and $O$ a compact normal subgroup such that $H/O$ is a Lie group. Let $\{\overline{x}(t)\}_t$ be a $1$-parameter subgroup of $H/O$.
Then $\{\overline{x}(t)\}_t$ can be lifted  along the projection
$H\rightarrow H/O$ to a $1$-parameter subgroup $\{x(t)\}_t$ in $H$.
\end{lem}

\begin{lem} \label{BGSS4.2} (\cite{BGSS} \textit{4.2})
Let $H$ be a locally compact group and $O$ a compact normal subgroup
such that $H/O$ is a connected Lie group.
Then $H=Z_H(O)O$.
\end{lem}

Using Lemma \ref{MZ4.6}, we can sharpen Lemma \ref{BGSS4.2} as follows:

\begin{cor}\label{cor:conn-centr}
Under the assumptions of Lemma \ref{BGSS4.2}, $H=Z_H(O)^\circ O$.
\end{cor} 

%We will also use:

%\begin{lem}\label{lem:nil-abel}
%If $H$ be a locally compact group admitting a compact central subgroup $C$ such that $H/C$ is an abelian Lie group, then $H$ is abelian.
%\end{lem}

We will also make use of the following classical:

\begin{thm}(\cite[page 137]{kaplansky})\label{thm:137}
Let $G$ be a connected locally compact group. Then $G$ admits a descending sequence of compact normal subgroups $N_i$ with trivial intersection such that $G/N_i$ is a Lie group for every $i$.
\end{thm}

Note that if $K\lhd G$ is any compact normal co-Lie subgroup, one may replace $N_i$ by $K\cap N_i$ and get such a sequence inside $K$. The normality and compactness of $K\cap N_i$ is obvious. Furthermore, since $K\cap N_i$ is compact 
$G/K\cap N_i$ is isomorphic to a {\it closed} connected subgroup of the Lie group $G/K\times G/N_i$, hence is also a Lie group.

\begin{proof}[Proof of Proposition \ref{prop:non-sol-quo}]
%Let $\gC$ be a limit group and consider the space $\text{Hom}(\gC,G)$ of representations of $\gC$ in $G$. Since $G$ is locally compact, $\text{Hom}(\gC,G)$ carries a complete uniform structure. Let $R$ be the solvable radical of $G/K$ and let $\ti R$ be its pre-image in $G$ under the quotient map. 
In view of Corollary \ref{cor:conn-centr} we may replace $H$ by $Z_H^\circ(O)$ and assume that $H$ is connected and $O$ is central.
Let $S$ be a semisimple Levi factor of the connected Lie group $H/O$, then $H/O=S\ltimes R^\circ$ is an almost semidirect product, and by the assumption $S$ is nontrivial. %Denote by $\ti R$ the pre-image of $R$ in $H$ under the quotient map. Then $\ti R$ is amenable. In fact $\ti R$ is the maximal abstractly amenable normal subgroup, i.e. every normal subgroup which is amenable as an abstract group is contained in $\ti R$, and $G/\ti R=\Ad(S)$. 

%Given a nonabelian limit group $\gC$ we have to produce an imbedding $H$ which projects faithfully and densely to the Adjoint group $\Ad(S)$.} 

%\medskip

%First note that $\Ad(S)$ admits a dense $F_2$ (cf. \cite{dense}), which, being free, can be lifted arbitrarily to $H$. 
%The claim is proved along the same strategy that we used in the previous sections. That is, given a limit subgroup $L\le H$ which projects densely and faithfully to $\Ad(S)$, and a generalized double $L'=A*_{C}B\xrightarrow{\varphi}L$ (or $L'=A*_{C}\xrightarrow{\varphi}L$) over $L$, we show that $L'$ can be imbedded in $H$ with trivial intersection with $\ti R$ and dense projection to $\Ad (S)$. 

By Theorem \ref{thm:137} $O$ admits a descending chain of compact subgroups $O_i$ with trivial intersection, such that $H/O_i$ is a Lie group for every $i$. Set $O_0=O$, and denote by $\pi_i:H\to H/O_i,~i\ge 0$ the corresponding projections, as well as, abusing notations, the projections $H/O_j\to H/O_i$ when $O_j\le O_i$.
Since, for $j>i$, the Lie algebra $\text{Lie}(H/O_j)$ is a central extension of $\text{Lie}(H/O_i)$, we can chose the Levi subgroups $S_i\le H/O_i$ so that $S_j$ covers $S_i$ whenever $j>i$. Let $\ti S$ be the universal covering group of $S$, and denote by $f_i:\ti S\to S_i$ the covering maps. Note that then for $j>i$, $f_i=\pi_i\circ f_j$.

Let $\gC$ be a nonabelian limit group. By Lemma \ref{lem:S->Ad(S)}, $\ti S$ admits a dense copy of $\gC$, which we will also denote by $\gC$, that projects faithfully to $\Ad_S(S)=\Ad_{\ti S}(\ti S)$. For $\gc\in \gC$ define 
$$
 \phi(\gc)=\cap_{i=0}^\infty \pi_i^{-1}(f_i(\gc)).
$$ 
Note that for $\gc\in\gC$ and $j>i$, 
$$
 \pi_j^{-1}(f_j(\gc))\subset \pi_j^{-1}(f_j(\gc))O_i=\pi_i^{-1}(f_i(\gc)),
$$ 
and since $\cap_{i=0}^\infty O_i=\{ 1\}$, $\phi(\gc)$ is a singleton. Since for each $i$, $\pi_i^{-1}(f_i(\gC))$ is a group, also the intersection $\cap_{i=0}^\infty \pi_i^{-1}(f_i(\gC))$ is a group. Moreover, since $\pi_0\circ \phi=f_0|_\gC$ and $f_0|_\gC$ is injective, we deduce that $\phi:\gC\to H$ is injective as well as a homomorphism, indeed
$$
 \phi=(\pi_0|_{\phi(\gC)})^{-1}\circ f_0|_\gC.
$$
By construction, $\phi(\gC)$ projects densely and faithfully to $\Ad(S)\cong H/\ti R$.
\end{proof}

%Clearly $\{\overline x(t)\}_t$ can be lifted uniquely to a one parameter subgroup $\{\overline x_i(t)\}_t$ of $S_i$, which by Lemma \ref{lem:lift-1-parameter} can be lifted to a one parameter subgroup $\{x_i(t)\}_t$ in $H$. Observe that if $j>i$ the images of $x_i(t)$ and $x_j(t)$ coincide in $G/O_i$, hence $x_i(t)x_j(t)^{-1}\in O_i$. It follows that the one parameter groups $\{x_i(t)\}_t$ converge uniformly in $H$. Let $\{x(t)\}_t$ be the limit one parameter subgroup of $H$.  

%Let $K$ be an identity neighborhood in $\ti Z$. Then $K$ is a Baire space. Additionally if $K$ is sufficiently small,the homomorphism $h_k:L'\to H$ projects densely to $\Ad(S)$ (see \cite{zuk} or \cite{dense}). Finally, for every $\gc\in L'\setminus\{1\}$ the set $\{k: h_k(\gc)\in\ti R\}$ is closed with empty interior. As in the previous sections, this implies that for some $k\in K$ the map $h_k$ is also faithful (even into $\Ad (S)$), and the claim is proved.

\medskip

We are now able to complete the proof of Theorem \ref{thm:nondiscrete}. First recall another classical fact about locally compact groups:

\begin{thm}(\cite{MZ}, 4.5)\label{MZ4.5}
Let $G$ be a locally compact group. There exists an open subgroup $H\leq G$ and a compact subgroup $O$, normal in $H$, such that $H/O$ is a connected Lie group.
\end{thm}

Suppose now that $G$ is a locally compact group containing a nondiscrete nonabelian free subgroup $F$. Replacing $G$ by the closure of $F$, we may assume $F$ is dense, and replacing $G$ further by an open subgroup, we may assume that $G/O$ is a connected Lie group for some compact normal subgroup $O\lhd G$. 
If $G/O$ is solvable, then $F\cap O$, being a nontrivial normal subgroup of $F$, is a nonabelian free group, and as $O$ is compact, it is nondiscrete. Hence Theorem \ref{thm:nondiscrete} follows, in this case, from Theorem \ref{thm:compact}. Suppose therefore that $G/O$ is nonsolvable. Then by Proposition \ref{prop:non-sol-quo}
there is a faithful homomorphism $f:\gC\to G$ such that $f(\gC)\ti R$ is dense in $G$, where $\ti R$ is the pre-image in $G$ of the solvable radical of $G/O$. But then, as $\ti R$ is amenable, and $G/\ti R$, being connected, has no nontrivial open subgroups, it follows from the generalized Auslander's theorem \cite[Theorem 1.13]{toti} that $f(\gC)$ is not discrete.
\qed

\begin{rem}
In the proof of Proposition \ref{prop:non-sol-quo}, one can actually produce a locally faithful homomorphism from $\ti S$ to $H$ whose image covers all the $S_i$'s, and whose closure is isomorphic to the inverse limit of the $S_i$'s. In this way one obtains directly a nondiscrete imbedding of $\gC$ in $H$ without using the generalized Auslander's theorem.
\end{rem} 

%%%%%%%%%%%%%%%%%%%%%%%%%%%%%%%%%%%%%%%%%%%%%%%%%

\section{Dense embeddings of surface groups}\label{sec:surface}

This section has two purposes. One is to correct a mistake from \cite{BGSS}.
The problem there is at the bottom of pg. 1380 where the map $\rho _{\beta}$ is ill defined; its image is not contained in $H/K$ but rather in the coset space $G/K$. The second is to prove Theorem \ref{thm:surface}, extending the main result from \cite{BGSS} to general oriented surfaces, not necessarily of even genus.

\subsection{Adding an independent element to a dense subgroup}\label{sec:infusion}
The proof of Theorem \ref{thm:surface} relies on the possibility of enlarging an independent set of group elements, assuming it generates a nondiscrete subgroup. The following proposition generalizes  Lemma \ref{infusion}. The proof of this fact given in \cite{BGSS} is inaccurate. 

\begin{prop} \label{prop:extension}
Let $G$ be a locally compact group and let ${x_1,\ldots,x_k}\subseteq G$ be
independent elements which generate a (rank $k$ free) group $F$ which
is not discrete in $G$. Let $U\subseteq G$ be an identity neighborhood.
Then there exists $x_{k+1}\in U$ such that $x_1,\ldots,x_{k+1}$ are
still independent.
\end{prop}

For the proof of Theorem \ref{thm:surface} we shall need to choose  $x_{k+1}$ more carefully; the possibility to do so will follow from the proof and is summarized in Corollary \ref{cor:extension}.

\begin{proof}
Let $F'$ be an abstract free group of rank  $k+1$, freely generated
by $\{y_1,\ldots,y_{k+1}\}$. We consider the homomorphisms $F'\rightarrow G$ which send $y_i\mapsto x_i$ for $i\leq k$ and $y_{k+1}\mapsto g$ for various choices of $g\in G$ and show that we may find a $g$ which yields an embedding. 
Remember that the identity component $G^\circ$
is closed, hence a locally compact, normal subgroup of $G$. We treat
the following three cases separately:

\medskip

\noindent\textbf{Case I} $F\cap G^\circ=\{ 1\}$. 
First assume $G^\circ=\{ 1\}$, i.e. that $G$ is totally disconnected, and let $K\leq G$ be a compact open subgroup. 
Since $F$ is nondiscrete, it intersects $K$ nontrivially, and hence admits a nondiscrete cyclic subgroup, so the result follows in this case from Lemma \ref{infusion}.

When $G$ is not totally disconnected $F$ injects (since $F\cap
G^\circ=\{1\}$) via the natural quotient into the totally disconnected
group $G/G^\circ$ and as above we may find a faithful map
$F'\rightarrow G/G^\circ$. Any lifting of the
image of $y_{k+1}$ from $G/G^\circ$ to $G$ yields a faithful embedding $F'\to G$ and completes the proof of the proposition in case I. Note that in this case although $x_{k+1}$ does not lie in a compact subgroup of $G$ it can be chosen to map to a compact subgroup in $G/G^\circ$. 
\medskip

In the next two cases we assume $F\cap G^\circ\neq \{1\}$, whence it
follows that $F$ intersects every identity neighborhood of $G^\circ$:
since $F$ is not discrete in $G$ and $G^\circ$ is normal, for $f\in F\cap G^\circ$ we may
choose $g\in F$ so that $[g,f]$ is nontrivial and arbitrarily close to the
identity. 

It follows from Theorem \ref{MZ4.5} that there exists a compact normal subgroups
$N\vartriangleleft G^\circ$ such that $G^\circ/N$ is a Lie group.

\medskip

\noindent\textbf{Case II} 
There is a compact subgroup $N\vartriangleleft G^\circ$, co-Lie, such that 
$F\cap N \neq \{1\}$: as in the first part of Case I we may find $x_{k+1}\in N$ arbitrarily close
to the identity.

\medskip

\noindent\textbf{Case III} $F\cap G^\circ\neq \{1\}$ but $F$ trivially intersects every compact normal co-Lie subgroup $N\vartriangleleft G^\circ$. 
Since $G^\circ$ is normal in $G$, for each $f\in F$ the conjugation by
$f$ is an automorphism of $G^\circ$, hence $G^\circ/N^f$ is also a Lie group.

Given a nontrivial reduced word $w$ in $k+1$ letters, we aim to show that the set
$$
 \mathcal {O}_w:= \{g\in G^\circ: w(x_1,\ldots,x_k,g)\neq 1\}
$$ 

\noindent 
is open dense in $G^\circ$. Let $\ti w$ be the
(continuous) map $\ti w:G^\circ\rightarrow G$ defined by
$$ 
 \ti w(g)=w(x_1,\ldots,x_k,g).
$$
Breaking up $w$ appropriately there exist words $w_j=w_j(y_1,\ldots,y_k)$
and $l_j\in \mathbb{Z}\setminus\{ 0\}$ such that
$$
 w(y_1,\ldots,y_k,y_{k+1})=w_my_{k+1}^{l_m}w_{m-1}\cdots w_1y_{k+1}^{l_1}w_0.
$$
For $i\le m$, denote
$$
 h_i=w_i(x_1,\ldots,x_k)w_{i-1}(x_1,\ldots,x_k)\cdots w_0(x_1,\ldots,x_k).
$$
Then
$$
 \ti w(g)=\ti w(1)\Pi_{i=1}^{m} (g^{l_{m-i+1}})^{h_{m-i}^{-1}}.
$$
Since $G^\circ\vartriangleleft G$, for $g\in G^\circ$ we have $\prod_{i=1}^{m}
(g^{l_{m-i+1}})^{h_{m-i}^{-1}}\in G^\circ$. Hence
$$
 \forall g\in G^\circ,~~\ti w(1)\in G^\circ\iff \ti w(g)\in G^\circ.
$$
If $\ti w(1)\notin G^\circ$ then $\forall g\in G$, $\ti w(g)\notin G^\circ$ 
and in particular $\mathcal{O}_w=G^\circ$ which is open dense in $G^\circ$.

\medskip
\noindent Assume $\ti w(1)\in G^\circ$. Then $\ti w:G^\circ\rightarrow G^\circ$ is
well defined. Fix a compact normal co-Lie subgroup $N\lhd G^\circ$.

\begin{lem}\label{word function}
There exists a compact normal subgroup $N_w\vartriangleleft G^\circ$ such that:
\begin{itemize}
\item $G^\circ/N_w$ is a Lie group. \item For $g,h\in G^\circ~~$  $gN_w=hN_w\implies w(g)N=w(h)N$. \item The induced map $\overline{w}:G^\circ/N_w\rightarrow
G^\circ/N$ is analytic.
\end{itemize}
\end{lem}

\begin{proof}
Set $N_w=\bigcap_{i=0}^m N^{h_i}$. Clearly $N_w$ is a compact normal
subgroup of $G^\circ$. Consider the map $G^\circ\to \prod_{i=0}^m G^\circ/ N^{h_i},~g\mapsto (gN^{h_i})_{i=1}^m$. Since the kernel $N_w$ is compact the image is closed. Because the target group $\prod_{i=0}^m G^\circ/ N^{h_i}$ is a Lie group, it follows that $G^\circ/N_w$ is a Lie group (recall that a topological group is a Lie group iff it is isomorphic to a Lie group in which case the analytic structure is uniquely defined, see \cite[Lemma 4.7.1]{MZ} or \cite[Theorem 5.3.2]{conlon}).

Let $n\in N_w$ and $g\in G^\circ$. Since $N_w$ is normal in $G^\circ$, for $i=0,\ldots,m$ there exists $n'_i\in N_w$ such that
$(gn)^{l_i}=g^{l_i}n'_i$. As $N_w\subseteq N^{h_i}$  we have
$$
 \ti w(gn)N=\ti w(1)\prod_{i=m}^1 ((gn)^{l_i})^{h_{i-1}^{-1}}N=\ti w(1)\prod_{m}^1 \left((g^{l_i})^{h_{i-1}^{-1}}{n'}_i^{h_{i-1}^{-1}}N\right)=\ti w(1)\prod_{m}^1 (g^{l_i})^{h_{i-1}^{-1}}N=\ti w(g)N.
 $$

It remains to show that the induced map $\overline{w}$ is analytic:
For each $i$, the map $gN_w\mapsto g^{h_i^{-1}}N$ is a well-defined continuous homomorphism
between Lie groups, hence analytic. Therefore
\begin{eqnarray*}
\overline{w}(gN_w)=\ti w(g)N = \ti w(1)\prod_{i=m}^1 \big( g^{h_{i-1}^{-1}}N\big)^{l_i}
\end{eqnarray*}
is analytic being the multiplication of analytic functions. This completes the proof of Lemma \ref{word function}.
\end{proof}

%------------------------------------------------------------------------

We now complete the proof of Proposition \ref{prop:extension}. By assumption, there exists
$a\in F\cap G^\circ$. Pick $i$ such that  $a,a^{x_i}$ do
not commute. By Corollary \ref{cor:baum1} for $l\in \mathbb{Z}$ large
enough, $\ti w((a^{x_i})^{-l}a(a^{x_i})^l)\neq 1$. Since $F\cap N=1$, $F\cap G^{\circ}$
injects into $G^\circ/N$ so 
$\overline{w}((a^{x_i})^{-n}a(a^{x_i})^n)\neq1$. Since $\overline{w}:G^\circ/N_w\to G^\circ/N$ is
analytic it follows that $(G^\circ/N_w)\setminus \overline{w}^{-1}(1)$ is open
dense and so is its preimage in $G^\circ$ (the quotient map is open). Finally note that this pre-image is contained in $\mathcal{O}_w$. 

Thus we have shown that for every nontrivial word $w$ in $k+1$ variables, the set $\mathcal{O}_w$ is open dense in $G^\circ$, hence by Baire's category theorem we conclude that 
$$
 \bigcap_{w\in F'} \mathcal{O}_w
$$
is dense in $G^\circ$, and in particular nonempty. For every $g\in \bigcap_{w\in F'} \mathcal{O}_w$ the map $\rho_g:F_{k+1}\rightarrow G$ defined by
$\rho_g(y_i)=x_i$ for $1\leq i\leq k$ and $\rho_g(y_{k+1})=g$ is
faithful and setting $x_{k+1}:=g$ completes the proof of the proposition in this case.
Moreover since the intersection is nondiscrete we may choose
$x_{k+1}\in G^\circ$ arbitrarily close to the identity. 
\end{proof}

The following consequence of the proof will be useful to us:

\begin{cor}\label{cor:extension}
Suppose that $W$ is a word in $k+1$ letters such that $W(y_1,\ldots,y_k,1)=1$ for $y_1,\ldots,y_k$ independent. Then, given $G$ and $x_1,\ldots,x_k$ as in the statement of the proposition, we can chose $x_{k+1}$ independent to $\{x_1,\ldots,x_k\}$, such that:
\begin{itemize}
\item In Case $I$, $\langle x_1,\ldots,x_{k+1}\rangle\cap G^\circ=\{1\}$, and $W(x_1,\ldots,x_k,x_{k+1})\in H$, where $H$ is an open subgroup of $G$ which is compact modulo $G^\circ$.

\item In Case $II$, $W(x_1,\ldots,x_{k+1})$ belongs to the compact normal subgroup $N$ of $G^\circ$.

\item In Case $III$, $W(x_1,\ldots,x_{k+1})N$ is contained in a $1$-parameter subgroup of $G^\circ/N$.  
\end{itemize} 
\end{cor}

Modulo Proposition \ref{prop:extension} and Corollary \ref{cor:extension} the proof given in \cite[Theorem 1.2]{BGSS} for Theorem \ref{thm:surface} for surfaces of even genus is correct, or more precisely can easily be corrected arguing similarly to Paragraph \ref{subsec:surface-proof} below. 
We will therefore concentrate below on the case of odd genus surfaces.
Note that every orientable surface covers an orientable
surface of even genus, hence the existence of odd genus embedding follows from the even genus. Our main issue is dealing with the density question.

\subsection{An eventually faithful sequence from $\Gamma_{2r+1}\rightarrow F_{2r+1}$}\label{sec:efs}

Let $\Gamma=\Gamma_{2r+1}$ be the fundamental group of the closed
orientable surface of genus $2r+1$ ($r\geq1$). We have the
presentation
\[
 \Gamma=\langle a_{1},a_{1}',\dots,a_{r},a_{r}',b,b',c_{1},c_{1}',\dots,c_{r},c_{r}'\mid[a_{1},a_{1}']\cdots[a_{r},a_{r}'][b',b]
 [{c}_{1}',c_{1}]\cdots[c_  {r}',c_{r}]\rangle.
\]

We shall produce a specific eventually faithful sequence of homomorphisms from $\Gamma$ to a fixed free group to be used in the proof of Theorem \ref{thm:surface}.

Denote 
$$
 \alpha=[a_{1},a_{1}']\cdots[a_{r},a_{r}']b'~\text{and}~
 \beta=b'~~\text{(cf. Figure \ref{surface})}.
$$

\begin{figure}[h!]

\begin{center}
\input{surface.pdf_t}
\caption{}
\label{surface}
\end{center}
\end{figure}

Define automorphisms $\sigma,\tau:\Gamma\rightarrow\Gamma$ to be the
Dehn twists (see Definition \ref{dehn twist}) around the boldface loops in Figure \ref{surface}. These can be thought of in terms of the effect of the topological twist on the fundamental group. Or algebraically, as the Dehn twist around $\alpha$ and $\beta$  with respect to appropriate HNN splittings of $\Gamma$. 
We give an explicit description as well, 
$$
\begin{matrix}
 \sigma(a_i)=a_i &  \tau(a_i)=a_i \\
 \sigma(a_i')=a_i' & \tau(a_i')=a_i' \\
 \sigma(b')=b' & \tau(b')=b' \\
  \sigma(b)=\alpha b & \tau(b)=b\beta ^{-1} \\ 
 \sigma(c_i)= \alpha c_i\alpha ^{-1} & \tau(c_i)=c_i \\
 \sigma(c_i')=\alpha c_i'\alpha ^{-1} & \tau(c_i')=c_i'.
\end{matrix}
$$
Note that $\sigma$ and $\tau$ commute.

Let $F_{2r+1}$ be an abstract free group with free generators $\{
x_{1},x_{1}',\ldots,x_{r},x_{r}',y'\}$ and let $f:\Gamma\rightarrow
F_{2r+1}$ be the map defined by
$f(a_{i})=f(c_{i})=x_{i},~~f(a_{i}')=f(c_{i}')=x_{i}',~~ f(b)=1,~~
f(b')=y'$;
this map is induced by the topological folding map of the surface
across $\alpha$ and $\beta$, when viewing $F_{2r+1}$ as the fundamental group of the half surface with two boundary components.  

Denote 
$$
 \delta=\sigma\circ\tau,\,\,f_{n}=f\circ\delta^{n}~\text{and}~
 y:=f(\alpha)=[x_{1},x_{1}']\cdots[x_{r},x_{r}']y'.
$$

\begin{lem} \label{faithful sequence}
The sequence $f_{n}$ is eventually faithful.
\end{lem}

\begin{proof}
We use Lemma \ref{baumslag}. For any $x\in F_{2r+1}$ we have
$[xyx^{-1},y'],[xy'x^{-1},y]\neq 1$. To see this, note that both $y$ and $y'$ are
primitive and $y$ is not conjugate to $y'^{\pm 1}$ (because the image of $y$ in $F_{2r+1}/\langle \langle y' \rangle \rangle$ is not trivial).
It follows that $y$ and $y'$ fulfill the conditions required of the
$z_i$'s in Lemma \ref{baumslag} for any choice of $u_i\mbox{'s}\in F_{2r+1}$ which do not  belong to either of the cyclic groups $\langle y\rangle,\langle y'\rangle$.

We now proceed to decompose $f_n (g)$ into a product of $u_i$'s and $z_j$'s as in Lemma \ref{baumslag}. Let $g\in \Gamma\setminus\{ 1\}$. Write $g$
as a word in $a_{i},a_{i}',b,b',c_{i},c_{i}'$ (this is not unique)
and break it into subwords
$$
 g=b^{k_{l+1}}w'_{l}b^{k_{l}}\cdots b^{k_{1}}w'_{0}b^{k_{0}},\,\,\, k_{i}\in\mathbb{Z}
$$
where either $w'_{i}=w'_{i}(a_{j},a_{j}',b')$ (i.e. a word in
$a_{j},a_{j}',b'$) or $w'_{i}=w'_{i}(c_{j},c_{j}')$, and if
$k_{i}=0$ then $w'_{i-1}$ and $w'_{i}$ cannot be written as words of
the same type; i.e. for every $i$,
$w'_{i}$ is realizable by a loop contained in a half surface and no
two adjacent $w'_{i}$'s are realizable by loops in the same half
surface. We then obtain a decomposition for $\delta^n(g)$:
$$
 \delta^{n}(g)=\left(\alpha^{n}b\beta^{-n}\right)^{k_{l+1}}\alpha^{\gep_{l}n}w'_{l}
 \alpha^{{-\gep}_{l}n}\left(\alpha^{n}b\beta^{-n}\right)^{k_{l}}\cdots\left(\alpha^{n}b
 \beta^{-n}\right)^{k_{1}}\alpha^{\gep_{0}n}w'_{0}\alpha^{-\gep{}_{0}n}\left(\alpha^{n}b\beta^{-n}\right)^{k_{0}}
$$
where $\gep_{j}\in\{0,1\}$ depending on whether $w'_{j}$ is in the upper 
half surface or not. The intermediate words between the $w'_{i}$'s are of
the form
$$
 \alpha^{-\gep_{j}n}(\alpha^{n}b\beta^{-n})^{k_{j}}\alpha^{\gep_{j-1}n}
$$
and $k_{j}=0$ implies that $\gep_{j-1}$ and $\gep_{j}$ are different (since
$w'_{i}$ and $w'_{i-1}$ lie in opposite half surfaces). In either case

$$
 f\left(\alpha^{-\gep_{j-1}n}(\alpha^{n}b\beta^{-n})^{k_{j}}\alpha^{\gep_{j}n}\right)=
 y^{-\gep_{j-1}n}(y^{n}y'^{-n})^{k_{j}}y^{\gep_{j}n}\neq 1.
$$

Next, we reduce the decomposition of $\delta^n(g)$ by combining powers of 
$\alpha$ that arise from $w'_{i}(a_i,a'_i,b')$'s with adjacent ones arising from the twists; and similarly for $\beta$.  Mapping to $F_{2r+1}$ we obtain a decomposition
$$f_{n}(g)=s_{m+1}u_{m}\cdots u_{1}s_1u_0s_0$$
where: 
\begin{itemize}
\item 

$u_i$ is the image of one of the $w'_j$'s. It is a fixed word not depending on $n$, which is non trivial since $f|_{\text{half surface}}$ is
injective. These play the role of $u_i$'s in Lemma \ref{baumslag}.

\item 
$s_i$ is an alternating product of $y^{\pm n+c}$'s and $y'^{\pm n+c}$'s which is the image of the intermediate words (mentioned above) combined with adjacent fixed powers coming from the $w'_j$'s. For $n$ large enough these are non trivial.
\end{itemize}

As mentioned, $y$ and $y'$ play the role of the $z_j$'s in Lemma \ref{baumslag}. By virtue of our decomposition $u_{i}$ is a (fixed) power of $y$ only if it is adjacent to nontrivial powers of
$y'$ on both sides; and vice versa.   Thus $u_i\notin C(\langle z_j \rangle)$ whenever they are neighbors in $f_n(g)$ and the conditions of  Lemma \ref{baumslag} (cf. Corollary \ref{cor:bslg}) are satisfied, and we conclude that for large $n$ the element $f_n(g)\neq 1$ and thus $f_n$ is an eventually faithful sequence.

\end{proof}

%---------------------------------------------------------------------------------------------------

\subsection{The proof of Theorem \ref{thm:surface}}\label{subsec:surface-proof} 
Consider $\Gamma_{2r+1}$ with the presentation
$$
 \langle a_1,a'_1,\ldots,a_r,a'_r,b,b',c_1,c'_1,\ldots,c_1,c'_r \mid [a_1,a'_1]\cdots[a_r,a'_r][b',b][c_1',c_1]\cdots [c_r',c_r]\rangle
$$
as in Section \ref{sec:efs}.

Let $G$ be a locally compact group and $x_1,\ldots,x_r$ independent elements generating a dense subgroup.
By Proposition \ref{prop:extension} we may find $\{x'_1,\ldots,x'_r,y'\}\subseteq G$
which together with the $x_1,\ldots,x_r$ are independent. Denote
$y=[x_1,x'_1]\cdots [x_r,x'_r]y'$ and let $Z_{G\times
G}(y,y')$ be the centralizer of $(y,y')$ in $G\times G$. For each $(g,h)\in Z_{G\times
G}(y,y')$ we have

\begin{eqnarray*}
[x_1,x_1']\cdots [x_r,x_r'][y',gh^{-1}]g[x_1',x_1]\cdots [x_r',x_r]g^{-1}&=&\\
&=&ygh^{-1}y'^{-1}hg^{-1}gy'y^{-1}g^{-1}\\
&=&1.
\end{eqnarray*}
Hence we may define a homomorphism $\rho_{(g,h)}:\Gamma \rightarrow G$ by setting:
\begin{itemize}
\item $ \rho _{(g,h)}(a_i)=x_i$, $\rho _{(g,h)}(a'_i)=x'_i,~\forall 1\leq i\leq r$.
\item $\rho _{(g,h)}(b)=gh^{-1}$, $\rho _{(g,h)}(b')=y'$.
\item $\rho _{(g,h)}(c_i)=gx_ig^{-1}$, $\rho _{(g,h)}(c'_i)=gx'_ig^{-1},~\forall 1\leq i\leq r$
\end{itemize}

Since $\rho _{(y,y')}$ is simply the Dehn twist map $f\circ \delta$ of Section \ref{sec:efs},  by Lemma \ref{faithful sequence}, $\rho_{(y^n,y'^n)}$ is an eventually faithful sequence.

The space $Z_{G\times G}(y,y')$ parametrizes homomorphisms $\Gamma \rightarrow G$ whose image contains $F$. Moreover, it contains the eventually faithful sequence $\rho_{(y^n,y'^n)}$. We need to find a point corresponding to a faithful embedding.
We will treat the three cases specified in Section \ref{sec:infusion} (cf. Corollary \ref{cor:extension}) separately. 

\medskip

\noindent \textbf{Case I:} Here we assume that the free group 
$\langle x_1,\ldots,x_r\rangle$ injects to $G/G^\circ$. By Corollary \ref{cor:extension} we may choose the new independent elements $x_1',\ldots,x_r',y'$ so that $y,y'\in H$ for some open subgroup $H\le G$ which is compact modulo $G^\circ$.

First assume $G$ is totally disconnected so $H$ is compact, hence so is $\overline{\langle (y,y')\rangle}$, and we may argue, as in the proof of Theorem \ref{thm:compact}, that the set of elements in $\overline{\langle (y,y')\rangle}$ which do not induce a faithful embedding is of first Baire category.

When $G$ is not totally disconnected we cannot assume that $(y,y')$ lies in a compact group. However we can still argue as follows: by Theorem \ref{MZ4.5} we may assume that $H$ admits a compact normal subgroup $K\triangleleft H$ such that $H/K$ is a connected Lie group. 

\medskip

%\noindent Recall the following facts concerning the structure of locally compact groups.

%\begin{lem} (\cite{MZ} \textit{end of 4.6})\label{MZ4.6} Let $H$ be a locally compact group and $K$ a compact normal subgroup such that $H/K$ is a Lie group. Let $\{\overline{x}(t)\}_t$ be a $1$-parameter subgroup of $H/K$.
%Then $\{\overline{x}(t)\}_t$ can be lifted  along the projection
%$H\rightarrow H/K$ to a $1$-parameter subgroup $\{x(t)\}_t$ in $H$.
%\end{lem}

%Under these conditions, if in addition $H/K$ is connected it follows that $H=H^\circ K$.

We now consider $K \times K$ which is compact and normal in  $H\times H$ and note that the quotient is connected and Lie.
Thus by Corollary \ref{cor:conn-centr}
$$
 H\times H=(Z_{H\times H}(K \times K)) ^\circ K\times K.
$$
It follows that there exist $(y_z,y_z')\in (Z_{H\times H}(K \times K)) ^\circ \subseteq G^\circ \times G^\circ $ and $(y_k,y_k')\in K\times K$ such that 
$$
 (y,y')=(y_k,y_k')(y_z,y_z').
$$

Note that $[(y_k,y'_k ),(y,y')]=1$ so that $(y^n_k,y'^n_k )$ all yield homomorphisms. Since $ \langle x_1,...x_r,x_1',...,x_r',y'\rangle \cap G^\circ =\{1\}$ the sequence of homomorphisms corresponding to $(y^n,y'^n)$ is eventually faithful also mod $G^\circ$. Because  $(y,y')=(y_k,y_k') \mbox{mod} G^\circ \times G^\circ$ this sequence coincides, mod $G^\circ \times G^\circ$, with the one corresponding to $(y_k^n,y_k'^n)$. It follows that the latter sequence $\rho _{(y_k^n,y^n_k)}$ is eventually faithful. Since $(y_k,y_k')$ generates a group with compact closure, we finish as above.  

\noindent \textbf{Case II:} By Corollary \ref{cor:extension} we may choose
$x_1',\ldots,x_r',y'\in N$ for some $N\vartriangleleft G^\circ$ compact and normal. It then follows that $y$
belongs to the compact subgroup $N':=N^{x_1}NN^{x_2}N\cdots
N^{x_r}N\leq G^\circ$ and so $(y,y')$ belongs to the compact $N'\times
N$. We conclude as in Case I.

\noindent \textbf{Case III:} Here $N\vartriangleleft G^\circ$, $G^\circ/N$ is a Lie
group, $\langle x_1,\ldots,x_r,x_1',\ldots,x_r',y'\rangle$ intersects $N$
trivially and $\langle x_1',\ldots,x_r',y'\rangle\subseteq G^\circ$, hence
$y\in G^\circ$. By Corollary \ref{cor:extension} the elements $x_i'$ and $y'$ may be chosen close enough to the
identity so that the projection of $(y,y')$ into the Lie group
$(G^\circ\times G^\circ)/(N\times N)$ is contained in a $1$-parameter subgroup.

Consider $Z^\circ:=(Z_{G^\circ\times G^\circ}(y,y'))^\circ$. This is a closed
connected subgroup of $G^\circ\times G^\circ$. For every compact normal subgroup
$N'\vartriangleleft G^\circ\times G^\circ$ such that $G^\circ\times G^\circ/N'$ is a Lie
group, the image of $Z^\circ$ is a closed, hence Lie, subgroup isomorphic
to $Z^\circ/Z^\circ\cap N'$.

Let $\gamma \in \Gamma\setminus\{1\}$. We wish to show that
$\OO_{\gamma}:=\{(g,h)\in Z^\circ : \rho_{(g,h)}(\gamma)\neq 1\}$ is open
dense in $Z^\circ$. We may write $\gamma = w(a_i,a_i',b,b',c_i,c_i')$.
For $(g,h)\in Z^\circ$ define the evaluation map
$$
 \dot w(g,h):=w(x_i,x_i',gh^{-1},gx_ig^{-1},gh^{-1},gx_i'g^{-1})=\rho_{(g,h)}(\gamma)\in G.
$$
As in the proof of Proposition \ref{prop:extension}, 
$$
 \forall (g,h)\in Z^\circ,~~ \dot w(1,1)\in G^\circ \iff \dot w(g,h)\in G^\circ.
$$
Therefore if $\dot w(1,1)\notin G^\circ $ then $\mathcal{O}_\gamma =Z^\circ $ is open dense. Assume $\dot w(1,1)\in G^\circ$, so that  $w$  induces a continuous map $\dot w:Z^\circ\rightarrow G^\circ$.  As in Lemma \ref{word function} there exists $N_w \vartriangleleft G^\circ$ contained in $N$ such that $G^\circ/N_w$ is a Lie group and so that $\dot w$ factors through the canonical projections to an analytic map
$$
 \xymatrix{Z^\circ \ar[r]^{\dot w} \ar[d]^{\pi} & G^\circ \ar[d]^{\pi}\\
 \frac{Z^\circ}{(N_w \times N_w)\cap Z^\circ} \ar[r]^{\tilde{w}} &
 \frac{G^\circ}{N}}
$$

\noindent Let $\{ \overline{x}(t)\}_t$ be the $1$-parameter subgroup in $G^\circ\times G^\circ/N\times N$ with $\overline  x(1)= (y,y')N\times N$.

Since, by Lemma \ref{BGSS4.2} 
$$
 \frac{G^\circ\times G^\circ}{N\times N}\cong \frac{Z_{G^\circ\times
 G^\circ}(N\times N)}{(N\times N)\cap Z_{G^\circ\times G^\circ}(N\times N) }
$$
we may view $\{ \overline{x}(t)\}_t$ as a $1$-parameter subgroup
of the right-hand side. There exists a lifting of $\{
\overline{x}(t)\}_t$ to a $1$-parameter subgroup, $\{x_w(t)\} _t$, of the Lie group
$\frac{Z_{G^\circ\times G^\circ}(N\times N)}{Z_{G^\circ\times G^\circ}(N\times
N)\cap(N_w\times N_w)}$.

By Lemma \ref{MZ4.6} we may now proceed to lift $\{x_w(t)\} _t$ to a $1$-parameter
subgroup $\{x(t)\}_t$ of $Z_{G^\circ\times G^\circ}(N\times N)$.
Since $\overline{x}(t)=x(t)N\times N$ 
$$
 (y,y')N\times N=\overline{x}(1)= x(1)N\times N
$$ 
and hence $(y,y')=x(1)(n,n')$ for some
$(n,n')\in N\times N$. For all $t$, $x(t)$ commutes with $x(1)$ and
with $N\times N$, hence with $(y,y')$, and so $x(t)\subseteq Z^\circ$.

For $l\in \mathbb{Z}$ we have $x(l)=(y^l,y'^l)(h,h')$ for some
$(h,h')\in N_w\times N_w$, hence
$\rho_{x(l)}(\gamma)N=\rho_{(y^l,y'^l)}(\gamma)N$. Since
$\rho_{(y^l,y'^l)}$ is eventually faithful and its image does not
intersect $N$ it follows that for $\gamma \neq 1$ and for integral $l\gg 0$ we have
$$
 w(x(l))=\rho_{x(l)}(\gamma )N=\rho_{(y^l,y'^l)}N\neq 1. 
$$
Finally
$\tilde{w}$ is analytic and not constantly trivial hence
$$
 \frac{Z^\circ}{Z^\circ\cap N_w\times N_w}\setminus \tilde{w}^{-1}(1)
$$ 
is open dense and so is its preimage $\OO_{\gamma}\subseteq Z^\circ$. By
the Baire category theorem $\bigcap_{\gamma \in \Gamma} \OO_{\gamma}$ is nonempty, and any $(g_0,h_0)$ in this intersection gives
rise to a faithful map $\rho_{(g_0,h_0)}:\Gamma \rightarrow G$ whose
image contains $F$, hence is dense.
\qed

%%%%%%%%%%%%%%%%%%%%%%%%%%%%%%%%%%%%%%%%%

\end{document}

%% file: endAction.pdf_t
\begin{picture}(0,0)%
\includegraphics{endAction.pdf}%
\end{picture}%
\setlength{\unitlength}{3854sp}%
\begingroup\makeatletter\ifx\SetFigFont\undefined%
\gdef\SetFigFont#1#2#3#4#5{%
  \reset@font\fontsize{#1}{#2pt}%
  \fontfamily{#3}\fontseries{#4}\fontshape{#5}%
  \selectfont}%
\fi\endgroup%
\begin{picture}(5817,5482)(76,-8134)
\put( 91,-7036){\makebox(0,0)[lb]{\smash{{\SetFigFont{11}{13.2}{\familydefault}{\mddefault}{\updefault}{\color[rgb]{0,0,0}$z_{j}^{\epsilon '}$}%
}}}}
\put(5176,-7576){\makebox(0,0)[lb]{\smash{{\SetFigFont{11}{13.2}{\familydefault}{\mddefault}{\updefault}{\color[rgb]{0,0,0}$z_{i}^{\epsilon}$}%
}}}}
\put(5491,-3661){\makebox(0,0)[lb]{\smash{{\SetFigFont{11}{13.2}{\familydefault}{\mddefault}{\updefault}{\color[rgb]{0,0,0}$z_{j}^{-\epsilon '}$}%
}}}}
\put(4231,-6811){\makebox(0,0)[lb]{\smash{{\SetFigFont{11}{13.2}{\familydefault}{\mddefault}{\updefault}{\color[rgb]{0,0,0}$V_i^{\epsilon}$}%
}}}}
\put(4681,-4201){\makebox(0,0)[lb]{\smash{{\SetFigFont{11}{13.2}{\familydefault}{\mddefault}{\updefault}{\color[rgb]{0,0,0}$V_j^{-\epsilon '}$}%
}}}}
\put(406,-3166){\makebox(0,0)[lb]{\smash{{\SetFigFont{11}{13.2}{\familydefault}{\mddefault}{\updefault}{\color[rgb]{0,0,0}$u_n\cdot z_{j}^{\epsilon '}$}%
}}}}
\put(991,-4021){\makebox(0,0)[lb]{\smash{{\SetFigFont{11}{13.2}{\familydefault}{\mddefault}{\updefault}{\color[rgb]{0,0,0}$V_{j,n}^{\epsilon '}$}%
}}}}
\put(721,-6406){\makebox(0,0)[lb]{\smash{{\SetFigFont{11}{13.2}{\familydefault}{\mddefault}{\updefault}{\color[rgb]{0,0,0}$V_j^{\epsilon '}$}%
}}}}
\put(2386,-6901){\makebox(0,0)[lb]{\smash{{\SetFigFont{11}{13.2}{\familydefault}{\mddefault}{\updefault}{\color[rgb]{0,0,0}$z_j^{\epsilon 't}\cdot V^{\epsilon}_i$}%
}}}}
\put(2206,-5551){\makebox(0,0)[lb]{\smash{{\SetFigFont{11}{13.2}{\familydefault}{\mddefault}{\updefault}{\color[rgb]{0,0,0}$u_n\cdot V_j^{\epsilon '}$}%
}}}}
\end{picture}%

%% file: surface.pdf_t
\begin{picture}(0,0)%
\includegraphics{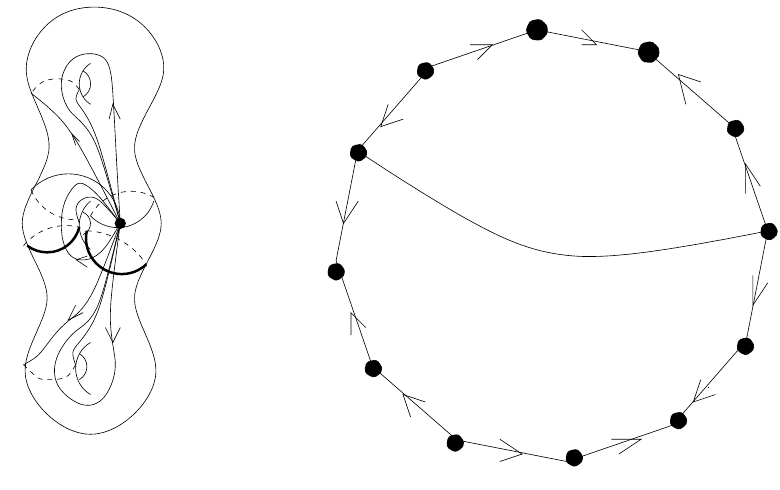}%
\end{picture}%
\setlength{\unitlength}{2611sp}%
\begingroup\makeatletter\ifx\SetFigFont\undefined%
\gdef\SetFigFont#1#2#3#4#5{%
  \reset@font\fontsize{#1}{#2pt}%
  \fontfamily{#3}\fontseries{#4}\fontshape{#5}%
  \selectfont}%
\fi\endgroup%
\begin{picture}(9409,5898)(706,-5530)
\put(721,-2626){\makebox(0,0)[lb]{\smash{{\SetFigFont{8}{9.6}{\familydefault}{\mddefault}{\updefault}{\color[rgb]{0,0,0}$\tau$}%
}}}}
\put(2566,-1996){\makebox(0,0)[lb]{\smash{{\SetFigFont{8}{9.6}{\familydefault}{\mddefault}{\updefault}{\color[rgb]{0,0,0}$\alpha$}%
}}}}
\put(9991,-1681){\makebox(0,0)[lb]{\smash{{\SetFigFont{8}{9.6}{\familydefault}{\mddefault}{\updefault}{\color[rgb]{0,0,0}$a$}%
}}}}
\put(9271,-511){\makebox(0,0)[lb]{\smash{{\SetFigFont{8}{9.6}{\familydefault}{\mddefault}{\updefault}{\color[rgb]{0,0,0}$a'$}%
}}}}
\put(7741,209){\makebox(0,0)[lb]{\smash{{\SetFigFont{8}{9.6}{\familydefault}{\mddefault}{\updefault}{\color[rgb]{0,0,0}$a$}%
}}}}
\put(6031, 29){\makebox(0,0)[lb]{\smash{{\SetFigFont{8}{9.6}{\familydefault}{\mddefault}{\updefault}{\color[rgb]{0,0,0}$a'$}%
}}}}
\put(5041,-871){\makebox(0,0)[lb]{\smash{{\SetFigFont{8}{9.6}{\familydefault}{\mddefault}{\updefault}{\color[rgb]{0,0,0}$b'$}%
}}}}
\put(4321,-2221){\makebox(0,0)[lb]{\smash{{\SetFigFont{8}{9.6}{\familydefault}{\mddefault}{\updefault}{\color[rgb]{0,0,0}$b$}%
}}}}
\put(4591,-3661){\makebox(0,0)[lb]{\smash{{\SetFigFont{8}{9.6}{\familydefault}{\mddefault}{\updefault}{\color[rgb]{0,0,0}$b'$}%
}}}}
\put(5311,-4831){\makebox(0,0)[lb]{\smash{{\SetFigFont{8}{9.6}{\familydefault}{\mddefault}{\updefault}{\color[rgb]{0,0,0}$b$}%
}}}}
\put(6661,-5461){\makebox(0,0)[lb]{\smash{{\SetFigFont{8}{9.6}{\familydefault}{\mddefault}{\updefault}{\color[rgb]{0,0,0}$c'$}%
}}}}
\put(8281,-5461){\makebox(0,0)[lb]{\smash{{\SetFigFont{8}{9.6}{\familydefault}{\mddefault}{\updefault}{\color[rgb]{0,0,0}$c$}%
}}}}
\put(9541,-4651){\makebox(0,0)[lb]{\smash{{\SetFigFont{8}{9.6}{\familydefault}{\mddefault}{\updefault}{\color[rgb]{0,0,0}$c'$}%
}}}}
\put(10081,-3301){\makebox(0,0)[lb]{\smash{{\SetFigFont{8}{9.6}{\familydefault}{\mddefault}{\updefault}{\color[rgb]{0,0,0}$c$}%
}}}}
\put(7741,-2491){\makebox(0,0)[lb]{\smash{{\SetFigFont{8}{9.6}{\familydefault}{\mddefault}{\updefault}{\color[rgb]{0,0,0}$\alpha$}%
}}}}
\put(1306,-1681){\makebox(0,0)[lb]{\smash{{\SetFigFont{8}{9.6}{\familydefault}{\mddefault}{\updefault}{\color[rgb]{0,0,0}$b'$}%
}}}}
\put(1486,-2941){\makebox(0,0)[lb]{\smash{{\SetFigFont{8}{9.6}{\familydefault}{\mddefault}{\updefault}{\color[rgb]{0,0,0}$b$}%
}}}}
\put(2161,-3976){\makebox(0,0)[lb]{\smash{{\SetFigFont{8}{9.6}{\familydefault}{\mddefault}{\updefault}{\color[rgb]{0,0,0}$c$}%
}}}}
\put(1306,-1141){\makebox(0,0)[lb]{\smash{{\SetFigFont{8}{9.6}{\familydefault}{\mddefault}{\updefault}{\color[rgb]{0,0,0}$a'$}%
}}}}
\put(2116,-466){\makebox(0,0)[lb]{\smash{{\SetFigFont{8}{9.6}{\familydefault}{\mddefault}{\updefault}{\color[rgb]{0,0,0}$a$}%
}}}}
\put(2521,-2941){\makebox(0,0)[lb]{\smash{{\SetFigFont{8}{9.6}{\familydefault}{\mddefault}{\updefault}{\color[rgb]{0,0,0}$\sigma$}%
}}}}
\put(766,-4066){\makebox(0,0)[lb]{\smash{{\SetFigFont{8}{9.6}{\familydefault}{\mddefault}{\updefault}{\color[rgb]{0,0,0}$c'$}%
}}}}
\end{picture}%